%
%
%

\def\LaTeX{%
  \let\Begin\begin
  \let\End\end
  \let\salta\relax
  \let\finqui\relax
  \let\futuro\relax}

\def\UK{\def\our{our}\let\sz s}
\def\USA{\def\our{or}\let\sz z}



\LaTeX

\USA


\salta

\documentclass[twoside,12pt]{article}
\setlength{\textheight}{24cm}
\setlength{\textwidth}{16cm}
\setlength{\oddsidemargin}{2mm}
\setlength{\evensidemargin}{2mm}
\setlength{\topmargin}{-15mm}
\parskip2mm


\usepackage{amsmath}
\usepackage{amsthm}
\usepackage{amssymb}
\usepackage[mathcal]{euscript}

\usepackage[usenames,dvipsnames]{color}
%
%
%
\def\gianni{\color{blue}}            
\def\pier{\color{red}}               
\def\danielle{\color{red}}           
\def\rev{\color{red}}
%
%
\let\gianni\relax
\let\pier\relax
\let\danielle\relax
\let\rev\relax


\bibliographystyle{plain}


%

\finqui

\def\Beq{\Begin{equation}}
\def\Eeq{\End{equation}}
\def\Bsist{\Begin{eqnarray}}
\def\Esist{\End{eqnarray}}

\def\Bthm{\Begin{theorem}}
\def\Ethm{\End{theorem}}

\def\Bprop{\Begin{proposition}}
\def\Eprop{\End{proposition}}
\def\Bcor{\Begin{corollary}}
\def\Ecor{\End{corollary}}
\def\Brem{\Begin{remark}\rm}
\def\Erem{\End{remark}}

\let\non\nonumber




\def\step #1 \par{\medskip\noindent{\bf #1.}\quad}

\def\Step #1 \par{\bigskip\leftline{\bf #1}\nobreak\medskip}


\def\Lip{Lip\-schitz}
\def\holder{H\"older}
\def\aand{\quad\hbox{and}\quad}

\def\lhs{left-hand side}
\def\rhs{right-hand side}
\def\sfw{straightforward}
\def\wk{well-known}


\def\characteriz{characteri\sz}

\def\organiz{organi\sz}

\def\regulariz{regulari\sz}

\def\bhv{behavi\our}


\def\multibold #1{\def\arg{#1}%
  \ifx\arg\pto \let\next\relax
  \else
  \def\next{\expandafter
    \def\csname #1#1#1\endcsname{{\bf #1}}%
    \multibold}%
  \fi \next}

\def\pto{.}

\def\multical #1{\def\arg{#1}%
  \ifx\arg\pto \let\next\relax
  \else
  \def\next{\expandafter
    \def\csname cal#1\endcsname{{\cal #1}}%
    \multical}%
  \fi \next}


\def\multimathop #1 {\def\arg{#1}%
  \ifx\arg\pto \let\next\relax
  \else
  \def\next{\expandafter
    \def\csname #1\endcsname{\mathop{\rm #1}\nolimits}%
    \multimathop}%
  \fi \next}

\multibold
qwertyuiopasdfghjklzxcvbnmQWERTYUIOPASDFGHJKLZXCVBNM.

\multical
QWERTYUIOPASDFGHJKLZXCVBNM.

\multimathop
dist div dom meas sign supp .


\def\accorpa #1#2{\eqref{#1}--\eqref{#2}}
\def\Accorpa #1#2 #3 {\gdef #1{\eqref{#2}--\eqref{#3}}%
  \wlog{}\wlog{\string #1 -> #2 - #3}\wlog{}}


\def\graffe #1{\mathopen\{#1\mathclose\}}

\def\<#1>{\mathopen\langle #1\mathclose\rangle}
\def\norma #1{\mathopen \| #1\mathclose \|}

\def\normaV #1{\norma{#1}_V}
\def\normaH #1{\norma{#1}_H}
\def\normaW #1{\norma{#1}_W}
\def\normaVp #1{\norma{#1}_*}

\def\iot {\int_0^t}
\def\ioT {\int_0^T}
\def\iO{\int_\Omega}
\def\intQt{\int_{Q_t}}
\def\intQ{\int_Q}

\def\dt{\partial_t}
\def\dn{\partial_\nu}

\def\cpto{\,\cdot\,}

\def\checkmmode #1{\relax\ifmmode\hbox{#1}\else{#1}\fi}
\def\aeO{\checkmmode{a.e.\ in~$\Omega$}}
\def\aeQ{\checkmmode{a.e.\ in~$Q$}}
\def\aet{\checkmmode{a.e.\ in~$(0,T)$}}
\def\aeS{\checkmmode{a.e.\ on~$\Sigma$}}
\def\aeG{\checkmmode{a.e.\ on~$\Gamma$}}

\def\aat{\checkmmode{for a.a.~$t\in(0,T)$}}


\def\erre{{\mathbb{R}}}




\def\genspazio #1#2#3#4#5{#1^{#2}(#5,#4;#3)}
\def\spazio #1#2#3{\genspazio {#1}{#2}{#3}T0}

\def\spaziom #1#2#3{\genspazio {#1}{#2}{#3}T{-1}}
\def\spaziomz #1#2#3{\genspazio {#1}{#2}{#3}0{-1}}
\def\spazioinf #1#2#3{\genspazio {#1}{#2}{#3}\infty0}
\def\L {\spazio L}
\def\H {\spazio H}
\def\W {\spazio W}

\def\LL {\spazioinf L}

\def\Lm {\spaziom L}
\def\Lmz {\spaziomz L}
\def\Hmz {\spaziomz H}
\def\C #1#2{C^{#1}([0,T];#2)}

\def\Cmz #1#2{C^{#1}([-1,0];#2)}

\def\Vp{V^*}
\def\Wz{W_0}
\def\Wzp{W_0^*}
\def\Hz{H_0}


\def\Lx #1{L^{#1}(\Omega)}
\def\Hx #1{H^{#1}(\Omega)}
\def\Wx #1{W^{#1}(\Omega)}

\def\Cx #1{C^{#1}(\overline\Omega)}

\def\Luno{\Lx 1}
\def\Ldue{\Lx 2}
\def\Linfty{\Lx\infty}
\def\Lq{\Lx4}
\def\Huno{\Hx 1}
\def\Hdue{\Hx 2}


\def\LQ #1{L^{#1}(Q)}


\let\a\alpha
\let\eps\varepsilon

\let\TeXchi\chi                         
\newbox\chibox
\setbox0 \hbox{\mathsurround0pt $\TeXchi$}
\setbox\chibox \hbox{\raise\dp0 \box 0 }
\def\chi{\copy\chibox}


\def\cd{c_\delta}
\def\ca{c_\a}
\def\cT{c_T}
\def\caT{c_{\a,T}}


\def\muz{\mu_0}

\def\uz{u_0}
\def\sigmaz{\sigma_0}
\def\mueps{\mu_\eps}
\def\ueps{u_\eps}
\def\sigmaeps{\sigma_\eps}
\def\Reps{R_\eps}
\def\muo{\mu_\omega}
\def\uo{u_\omega}
\def\sigmao{\sigma_\omega}
\def\mui{\mu_\infty}
\def\ui{u_\infty}
\def\sigmai{\sigma_\infty}
\def\xii{\xi_\infty}
\def\mus{\mu_s}
\def\us{u_s}
\def\sigmas{\sigma_s}
\def\xis{\xi_s}
\def\mua{\mu_\a}
\def\ua{u_\a}
\def\sigmaa{\sigma_\a}
\def\xia{\xi_\a}
\def\Ra{R_\a}
\def\mun{\mu_n}
\def\un{u_n}
\def\sigman{\sigma_n}
\def\xin{\xi_n}
\def\Rn{R_n}
\def\vn{v_n}
\def\tn{t_n}
\def\vstar{v^*}
\def\zstar{z^*}
\def\wstar{w^*}
\def\dh{\delta_h}
\def\vh{\bar v_h}
\def\muh{\bar\mu_h}

\def\Rh{\bar R_h}

\def\Beta{\widehat{\vphantom t\smash\beta\mskip2mu}\mskip-1mu}

\def\Betaeps{\Beta_\eps}
\def\betaeps{\beta_\eps}
\def\Weps{\calW_\eps}

\def\todx{\mathrel{\scriptstyle\searrow}}
\def\tosn{\mathrel{\scriptstyle\nearrow}}

\Begin{document}

\title{\bf On a Cahn-Hilliard type phase field system
related to tumor growth\footnote{{\bf Acknowledgments.}
The first two authors gratefully acknowledge some financial support from the MIUR-PRIN Grant 2010A2TFX2 ``Calculus of variations''.}}
\author{}
\date{}
\maketitle
\begin{center}
\vskip-2cm
{\large\bf Pierluigi Colli$^{(1)}$}\\
{\normalsize e-mail: {\tt pierluigi.colli@unipv.it}}\\[.25cm]
{\large\bf Gianni Gilardi$^{(1)}$}\\
{\normalsize e-mail: {\tt gianni.gilardi@unipv.it}}\\[.25cm]
{\large\bf Danielle Hilhorst$^{(2)}$}\\
{\normalsize e-mail: {\tt Danielle.Hilhorst@math.u-psud.fr}}\\[.45cm]
$^{(1)}$
{\small Dipartimento di Matematica ``F. Casorati'', Universit\`a di Pavia}\\
{\small via Ferrata 1, 27100 Pavia, Italy}\\
$^{(2)}$
{\small Laboratoire de Math\'ematiques}\\
{\small Universit\'e de Paris-Sud, 91405 Orsay, France}\\[.8cm]
\end{center}


\Begin{abstract}
The paper deals with a phase field system of Cahn-Hilliard type. For positive viscosity coefficients, the authors prove {\rev an existence and uniqueness} 
result and study the long time \bhv\ of the solution by assuming the nonlinearities to be rather general. In a more restricted setting, the limit as the viscosity coefficients tend to zero is investigated as well.\\[2mm]
{\bf Key words:} phase field model, {\danielle tumor growth}, viscous Cahn-Hilliard equations, well posedness, long-time \bhv, asymptotic analysis\\[2mm]
{\bf AMS (MOS) Subject Classification:} 35B20, 35K20, 35K35, 35R35
\End{abstract}


\salta

\pagestyle{myheadings}
\newcommand\testopari{\sc \ Colli \ --- \ Gilardi \ --- \ Hilhorst}
\newcommand\testodispari{\sc  Phase field system
related to tumor growth}
\markboth{\testodispari}{\testopari}

\finqui


\section{Introduction}
\label{Intro}
\setcounter{equation}{0}
In this article, we study
{\pier the coupled system of partial differential equations}
\Bsist
  && \a \dt\mu + \dt u - \Delta\mu
  = p(u) (\sigma - \gamma\mu)
  \label{Iprima}
  \\
  && \mu = \a\dt u - \Delta u + \calW'(u)
  \label{Iseconda}
  \\
  && \dt\sigma - \Delta\sigma
  = - p(u) (\sigma - \gamma\mu)
  \label{Iterza}
\Esist
{\pier in a domain $\Omega \times (0,\infty)$, together with the boundary conditions}
\Beq
  \dn\mu = \dn u = \dn\sigma = 0 {\pier \quad\hbox{on the boundary}}~\Gamma \times (0,\infty)
  \label{Ibc}
  \Eeq
{\pier and the initial conditions}
\Beq
  \mu(0) = \muz, \quad
  u(0) = \uz
  \aand
  \sigma(0) = \sigmaz \,.
  \label{Icauchy}
\Eeq
\Accorpa\Ipbl Iprima Icauchy
Each of the partial differential equations \accorpa{Iprima}{Iterza}
is meant to hold in a three-dimensional bounded domain $\Omega$
endowed with a smooth boundary~$\Gamma$ and for every positive time,
and $\dn$ in \eqref{Ibc} stands for the normal derivative.
Moreover, $\a$~and $\gamma$ are positive constants.
Furthermore, $p$~is a nonnegative function and {\danielle ${\cal W}$} is a nonnegative double well potential.
Finally, $\muz$, $\uz$, and $\sigmaz$ are given initial data defined in~$\Omega$.

The physical context is that of a tumor-growth model which has been derived from
{\rev a} general mixture theory {\rev \cite{OdeHawPruM3AS2010, HDVDZO}}. We also refer to
\cite{CriFriLiLowMacSanWisZheBOOK-CH2008}, \cite{LowFriJinChuLiMacWisCriNL2010} and \cite{CriLowBOOK2010} for overview articles and to \cite{CriLowNieJMB2003} and
\cite{CriLiLowWisJMB2009} for the study of related sharp interface models.

{\rev We point out that the unknown function $u$ is an order parameter which is close to two values in the regions of nearly pure phases, say $u \simeq 1$ in the tumorous phase and $u \simeq -1$ in the healthy cell phase; the second unknown 
$\mu$ is the related chemical potential, specified by \eqref{Iseconda} as in the case of the viscous Cahn-Hilliard or Cahn-Hilliard equation depending on whether $\alpha>0 $ or $\alpha =0$ (see \cite{CH, EllSt, EllZh});  
the third unknown $\sigma$ stands for the nutrient concentration, typically $\sigma \simeq 1$ in a nutrient-rich extracellular water phase and $\sigma \simeq 0$ in a nutrient-poor extracellular water phase.}  

{\rev In the case that the parameter $\alpha$ is strictly positive, {\danielle the problem~(\ref{Iprima})--(\ref{Icauchy})} is a generalized phase field model, while it becomes of pure Cahn-Hilliard type in the case that $\alpha = 0$. On the other hand, the presence of the term $\alpha \mu_t$ in \eqref{Iprima} gives, in the case $\alpha>0$, a parabolic structure to equation~\eqref{Iprima} with respect to $\mu$. Let us note that the meaning of the coefficient $\alpha$ here differs from the one in \eqref{Iseconda}: in \eqref{Iprima} $\alpha$ is not a viscosity coefficient since it enters in the natural Lyapunov functional of the system, which 
reads~(cf.~\cite{HKNZ})
$$
  E(u,\mu,\sigma ) = \int_\Omega \left( \frac12 |\nabla u|^2 + \calW (u) + \frac{\alpha}2 \mu^2 +\frac12 \sigma^2 \right). 
$$
However, the fact that the coefficients are taken equal in the system \accorpa{Iprima}{Iterza} is somehow related to the limiting problem obtained by formal asymptotics on $\alpha$. Indeed, {\danielle we refer to the forthcoming article \cite{HKNZ} for a formal study of the relation between these models and the corresponding sharp interfaces limits}.}

We remark that the original model deals with functions
$\cal W$~and $p$ that are precisely related to each other. Namely, we~have
\Beq
  p(u) = 2p_0 \sqrt{\calW(u)}
  \quad \hbox{if $|u|\leq 1$}
  \aand
  p(u) = 0
  \quad \hbox{otherwise}
  \label{relazWp}
\Eeq
where $p_0$ is a positive constant and where $ {\cal W}(u):=-\int_0^u f(s)\,ds$ is the classical
Cahn--Hilliard double well free-energy density. However, this relation is not used in this paper{\danielle , whose first aim is proving the well-posedness
of the initial-boundary value problem~(\ref{Iprima})--(\ref{Icauchy}) in the case of a positive parameter $\a$. In this setting, we can allow $\calW$ to be even a singular potential (the reader can see the later Remark~\ref{CompatWp}).

Actually, we prove the existence of a unique strong solution to the system~(\ref{Iprima})--(\ref{Icauchy}) under very general conditions on $p$ and ${\cal W}$, as well as we study the long time \bhv\ of the solution;
in particular, we can characterize the {\rev omega limit} set and deduce an interesting property in a special physical case (cf. the later Corollary~\ref{Specialcase}). Next, in a more restricted setting for the double-well potential ${\cal W}$, we investigate the asymptotic \bhv\ of the problem as the coefficient $\a$ tends to zero and find the convergence of subsequences to weak solutions of the limiting problem. Moreover, under a smoothness condition on the initial value $u_0$ we are able to show uniqueness for the limit problem and consequently also the convergence of the entire family as $\a \searrow 0$ (see Theorem~\ref{Asymptotics}).}

Our paper is \organiz ed as follows. In Section 2, we state our assumptions and results on the
mathematical problem. The forthcoming sections are devoted to the corresponding proofs. In Section
3, we prove the uniqueness of the solution. After presenting a priori estimates in Section 4, we prove the existence
of the solution on an arbitrary time interval in Section 5, while we study its large time behavior in Section 6;
Section 7 is devoted to the study of the limit of the phase field model (1.1)-(1.5) to the corresponding Cahn-Hilliard
problem as $\alpha \to 0$.


\section{Statement of the mathematical problem}
\label{Statement}
\setcounter{equation}{0}

In this section, we make {\pier precise assumptions} and state our results.
First of all,
we assume $\Omega$ to be a bounded connected open set in~$\erre^3$
(lower-dimensional cases could be considered with minor changes)
whose boundary~$\Gamma$ is supposed to be smooth.
As in the Introduction, the symbol $\dn$ denotes the (say, outward) normal derivative on~$\Gamma$.
As the first aim of {\pier our analysis}  is {\danielle to study the} well-posedness on any finite time interval,
we fix a final time $T\in(0,+\infty)$
and {\gianni let}
\Beq
  Q := \Omega\times(0,T)
  {\gianni \aand}
  {\pier \Sigma := \Gamma \times (0,T)} .
  \label{defQ}
\Eeq
Moreover, we set for convenience
\Beq
  V := \Huno,
  \quad H := \Ldue
  \aand
  W := \graffe{v\in\Hdue:\ \dn v = 0 \ \hbox{on $\Gamma$}}
  \label{defspazi}
\Eeq
and endow the spaces \eqref{defspazi} with their standard norms,
for which we use a self-explanato\-ry notation like $\normaV\cpto$.
{\pier If $p\in[1,+\infty]$, it will be useful to}
write $\norma\cpto_p$ for the usual norm
in~$L^p(\Omega)$.
In the sequel, the same symbols are used for powers of the above spaces.
It is understood that $H\subset\Vp$ as usual, i.e.,
in order that $\<u,v>=\iO uv$ for every $u\in H$ and $v\in V$,
where $\<\cpto,\cpto>$ stands for the duality pairing between $\Vp$ and~$V$.

\bigskip

As far as the structure of the system is concerned, we are given
two constants $\a$ and~$\gamma$ and three functions $p$, $\Beta$ and $\lambda$
satisfying the conditions listed below
\Bsist
  && \a \in (0,1) \aand \gamma > 0
  \label{hpconst}
  \\
  && \hbox{$p:\erre\to\erre$ is nonnegative, bounded and \Lip\ continuous}
  \label{hpp}
  \\
  && \hbox{$\Beta:\erre\to[0,+\infty]$ is convex, proper, lower
      semicontinuous}
    \qquad
  \label{hpBeta}
  \\
  && \hbox{$\lambda\in C^1(\erre)$ is nonnegative and
      $\lambda'$ is \Lip\ continuous}.
  \label{hplambda}
\Esist
We define the potential $\calW\!:\erre\to[0,+\infty]$
and the graph $\beta$ in $\erre\times\erre$ by
\Beq
  \calW := \Beta + \lambda
  \aand
  \beta := \partial\Beta
  \label{defWbeta}
\Eeq
\Accorpa\HPstruttura hpconst defWbeta
and note that $\beta$ is maximal monotone.
In the sequel, we write $D(\Beta)$ and $D(\beta)$ for the effective domains
of~$\Beta$ and~$\beta$, respectively, and we use the same symbol $\beta$
for the maximal monotone operators induced on $L^2$ spaces.

\Brem
\label{CompatWp}
Note that lots of potentials $\calW$ fit our assumptions.
Typical examples are the classical double well potential
and the logarithmic potential defined~by
\Bsist
  && \calW_{cl}(r) := {\textstyle \frac 14} (r^2-1)^2
  = {\textstyle \frac 14} ((r^2-1)^+)^2 + {\textstyle \frac 14} ((1-r^2)^+)^2
  \quad \hbox{for $r\in\erre$}
  \label{clW}
  \\
  && \calW_{log}(r) := (1-r)\ln(1-r) + (1+r)\ln(1+r) + \kappa (1 - r^2 )^+
  \quad \hbox{for $|r|<1$}
  \qquad
  \label{logW}
\Esist
where the decomposition $\calW=\Beta+\lambda$ as in \eqref{defWbeta} is written explicitely.
More precisely, in \eqref{logW}, $\kappa$~is a positive constant
which does or does not provide a double well depending on its value,
and the definition of the logarithmic part of $\calW_{log}$
is extended by continuity at $\pm1$ and by $+\infty$ outside~$[-1,1]$.
Moreover, another possible choice is the following
\Beq
  \calW(r) := I(r) + ((1 - r^2)^+)^2
  \quad \hbox{where $I$ is the indicator function of $[-1,1]$}
  \label{irrW}
\Eeq
{\pier taking the value $0$ in $[-1,1] $ and $+\infty$ elsewhere.}
Clearly, if $\beta$ is multi-valued like in the case~\eqref{irrW},
the precise statement of problem \Ipbl\
has to introduce a selection $\xi$ of~$\beta(u)$.
We also remark that our assumptions do not include
the relationship \eqref{relazWp} between $\rev \calW$ and~$p$,
and one can wonder whether what we have required is compatible with~\eqref{relazWp}.
This is the case if $\Beta$ and $\lambda$ satisfy suitable conditions, in addition.
For instance, we can assume the following:
$D(\Beta)$ includes the interval $[-1,1]$ and $\Beta$ vanishes there;
$\lambda$~is strictly positive in $(-1,1)$ and $\lambda(\pm1)=\lambda'(\pm1)=0$.
In such a case, $\calW$~presents two minima with quadratic \bhv\ at $\pm1$
and the function $p$ given by \eqref{relazWp} actually satisfies~\eqref{hpp}.
We note that this excludes the case of the logarithmic potential,
while it includes both \eqref{clW} and~\eqref{irrW}.
\Erem

As far as the initial data of our problem are concerned,
we assume that
\Beq
  \muz , \uz , \sigmaz \in V
  \aand
  \Beta(\uz) \in \Luno
  \label{hpdati}
\Eeq
{\pier while the regularity} {\danielle properties which we obtain for the solution
are the following}
\Bsist
  && \mu , u, \sigma \in \H1H \cap \L2W \subset \C0V
  \label{regcomponenti}
  \\
  && \xi \in \L2H
  \aand
  \xi \in \beta(u)
  \quad \aeQ .
  \label{regxi}
\Esist
\Accorpa\Regsoluz regcomponenti regxi
At this point, the problem we want to investigate consists
in looking for a quadruplet $(\mu,u,\sigma,\xi)$
satisfying the above regularity and the following boundary value problem
\Bsist
  \hskip-1cm & \a \dt\mu + \dt u - \Delta\mu = R
  \quad \hbox{where} \quad
  R = p(u) (\sigma - \gamma\mu)
  & \quad \aeQ
  \label{prima}
  \\
  \hskip-1cm & \mu = \a\dt u - \Delta u + \xi + \lambda'(u)
  & \quad \aeQ
  \label{seconda}
  \\
  \hskip-1cm & \dt\sigma - \Delta\sigma
  = -R
  & \quad \aeQ
  \label{terza}
  \\
  \hskip-1cm & {\pier \dn\mu = \dn u = \dn\sigma = 0}
  & {\pier \quad \aeS }
  \label{bc}
  \\
  \hskip-1cm & \mu(0) = \muz, \quad
  u(0) = \uz
  \aand
  \sigma(0) = \sigmaz
  & \quad \hbox{in $\Omega$} \,.
  \label{cauchy}
\Esist
\Accorpa\Pbl prima cauchy
\Accorpa\Tuttopbl regcomponenti cauchy
We note once and for all that {\rev adding} \eqref{prima} and~\eqref{terza} {\danielle yields}
\Beq
  \dt \bigl( \a \mu + u + \sigma \bigr)
  - \Delta (\mu + \sigma)
  = 0
  \quad \aeQ .
  \label{unomenotre}
\Eeq
Here is our well-posedness result.

\Bthm
\label{Wellposedness}
Assume \HPstruttura\ and \eqref{hpdati}.
Then, there exists a unique quadruplet $(\mu,u,\sigma,\xi)$
satisfying \Regsoluz\ and solving problem \Pbl.
\Ethm

\Brem
\label{Furtherregularity}
By starting from the regularity requirements \Regsoluz\
and owing to the regularity theory of elliptic and parabolic equations,
one can easily derive further properties of the solution.
For instance, {\pier as $R\in \L {\infty}H$}
one can show that
\Beq
  \sigma \in \W{1,p}H \cap \L pW
  \quad \hbox{for every $p\in[1,+\infty)$}
  \label{regsigma}
\Eeq
provided that $\sigmaz$ is smooth enough {\gianni(see, e.g., \cite[Thm.~2.3]{DHP})}.
\Erem

Once we know that there exists a unique solution on any finite time interval,
we can address its long time \bhv.
Our next result deals with the {\rev omega limit} of an arbitrary initial datum satisfying~\eqref{hpdati}.
Even though the possible topologies are several,
by recalling \eqref{regcomponenti} we choose
\Beq
  \Phi := V \times V \times V
  \label{phasespace}
\Eeq
as a phase space and~set
\Bsist
  \omega = \omega(\muz,\uz,\sigmaz) :=
  \Bigl\{ (\muo,\uo,\sigmao) \in \Phi :\ (\mu,u,\sigma)(\tn)\to
  (\muo,\uo,\sigmao)\qquad
  \non
  \\
  \quad {\pier \hbox{strongly in $\Phi$ for some $\{\tn\}\tosn+\infty$} \Bigr\}.}
  \label{omegalimit}
\Esist
Our result {\pier reads} as follows

\Bthm
\label{Longtime}
Assume \HPstruttura\ and \eqref{hpdati}.
Then, the {\rev omega limit} {\danielle set} $\omega$ is non-empty.
Moreover, if $(\muo,\uo,\sigmao)$ is any element of~$\omega$,
then $\muo$ and $\sigmao$ are constant
and their constant values $\mus$, $\sigmas$ and the function $\uo$
are related to each other~by
\Bsist
   p(\uo) (\sigmas - \gamma\mus) = 0
  \aand
  -\Delta\uo + \beta(\uo) + \lambda'(\uo) \ni \mus
  \quad \aeO , \qquad
  \non
  \\
\quad
  {\pier \dn\uo = 0 \quad \aeG.}
  \label{longtime}
\Esist
\Ethm

We deduce an interesting consequence in a special case which, however, is significant.

\Bcor
\label{Specialcase}
Assume that $D(\Beta)=[-1,1]$ and that
$p$ is strictly positive in $(-1,1)$
and vanishes at~$\pm1$.
Then, we have
\Beq
  \hbox{either} \quad
  \sigmas = \gamma\mus
  \quad \hbox{or} \quad
  \uo = -1
  \quad \aeO
  \quad \hbox{or} \quad
  \uo = 1
  \quad \aeO \,.
  \non
\Eeq
\Ecor

Indeed, if $\sigmas\not=\gamma\mus$,
\eqref{longtime} implies $\uo\in D(\Beta)$ and $p(\uo)=0$ \aeO.
On the other hand, we have $\uo\in V$, whence the conclusion.
We also remark that, in the case of the potential~\eqref{irrW},
the {\pier inclusion in} \eqref{longtime} reduces to $\beta(\uo)\ni\mus$.
In particular, $\uo=-1$ if $\mus<0$ and $\uo=1$ if $\mus>0$.

\medskip

Our final result regards the asymptotic analysis
as the viscosity coefficient $\a$ tends to zero
and the study of the limit problem.
We can deal with this by restricting ourselves
to a particular class of potentials.
Namely, we also assume~that
\Bsist
  && D(\Beta) = \erre
  \aand
    {\gianni \hbox{$\calW=\Beta+\lambda$ \enskip is a $C^3$ function on $\erre$}}
  \label{ovunque}
  \\
  && \calW(r) \geq \delta_0|r| - c_0 \,, \quad
  |\calW'(r)| \leq c_1 (|r|^3 + 1) \,, \quad
  |\calW''(r)| \leq c_2 (|r|^2 + 1)
  \qquad
  \label{classico}
\Esist
for any $r\in\erre$ and some positive constants $\delta_0$ and $c_0\,,c_1\,,c_2$.
We have written both the {\pier last two conditions in} \eqref{classico} for convenience
even though the latter implies the former.
We also remark that the classical potential~\eqref{clW} fulfils such assumptions.
Here is our result.

\Bthm
\label{Asymptotics}
Assume \HPstruttura, \eqref{hpdati}, and \accorpa{ovunque}{classico} in addition.
Moreover, let $(\mua,\ua,\sigmaa,\xia)$ be the unique solution to problem \Pbl\
given by Theorem~\ref{Wellposedness}.
Then, we have:
$i)$~the following convergence holds
\Bsist
  & \mua \to \mu
  & \hbox{weakly in $\L2V$}
  \label{convmu}
  \\
  & \ua \to u
  & \hbox{weakly star in $\L\infty V\cap\L2W$}
  \label{convu}
  \\
  & \sigmaa \to \sigma
  & \hbox{weakly in $\H1H\cap\L2W$}
  \label{convsigma}
  \\
  & \dt (\a\mua+\ua) \to \dt u
  & \hbox{weakly in $\L2\Vp$}
  \label{convdtamuu}
\Esist
at least for a subsequence;
$ii)$~every limiting triplet $(\mu,u,\sigma)$
satisfies
\Bsist
  & \< \dt u, v > + \iO \nabla\mu \cdot \nabla v
  = \iO R \, v
  & \quad \hbox{{\pier $\forall\,v\in V$,  \ \aet}}
  \label{primalim}
  \\
  & R = p(u) (\sigma - \gamma\mu)
  & \quad \aeQ
  \label{defR}
  \\
  & \mu = - \Delta u + \calW'(u)
  & \quad \aeQ
  \label{secondalim}
  \\
  & \dt\sigma - \Delta\sigma
  = - R
  & \quad \aeQ
  \label{terzalim}
  \\
  & {\pier \dn u = \dn\sigma = 0}
  & {\pier \quad \aeS}
  \label{bclim}
  \\
  & u(0) = \uz
  \aand
  \sigma(0) = \sigmaz
  & \quad \hbox{in $\Omega$} \,;
  \label{cauchylim}
\Esist
$iii)$~if
\Beq
  {\pier \uz \in W }
  \label{hpdatibis}
\Eeq
then every solution to problem \accorpa{primalim}{cauchylim}
satisfying the regularity given by \accorpa{convmu}{convdtamuu}
also satisfies the further regularity
\Bsist
  && \mu \in \L\infty H \cap \L2W
  \label{regmubis}
  \\
  && u \in \H1H \cap \L\infty W
  \subset \LQ\infty \,.
  \label{regubis}
\Esist
Moreover, such a solution is unique.
\Ethm

\Brem
\label{RegulUniq}
We observe that even further regularity for $\sigma$ could be derived
on account of the regularity of $R$ induced by~\eqref{regmubis},
provided that $\sigmaz$ is smoother.
It must be pointed out that a uniqueness result for the limit problem
has been proved in~\cite{FGR} by a different argument.
In the same paper, {\pier a slightly different regularity result is shown as well}.
Finally, we remark that {\danielle the uniqueness property} implies that
the whole family $(\mua,\ua,\sigmaa)$
converges to $(\mu,u,\sigma)$ in the sense of~\accorpa{convmu}{convdtamuu}
as $\a\todx0$.
\Erem

The rest of the section is devoted to list some facts
and to fix some notations.
We recall that $\Omega$ is bounded and smooth.
So, throughout the paper,
we owe to some \wk\ embeddings of Sobolev type,
namely
$V\subset\Lx p$ for $p\in[1,6]$,
together with the related Sobolev inequality
\Beq
  \norma v_p \leq C \normaV v
  \quad \hbox{for every $v\in V$ and $1\leq p \leq 6$}
  \label{sobolev}
\Eeq
and $\Wx{1,p}\subset\Cx0$ for $p>3$, together with
\Beq
  \norma v_\infty \leq C_p \norma v_{\Wx{1,p}}
  \quad \hbox{for every $v\in\Wx{1,p}$ and $p>3$}.
  \label{sobolevbis}
\Eeq
In \eqref{sobolev}, $C$ {\danielle only depends} on~$\Omega$,
while $C_p$ in \eqref{sobolevbis} {\danielle also depends} on~$p$.
In particular, the continuous embedding
$W\subset\Wx{1,6}\subset\Cx0$ holds.
Some of the previous embeddings are in fact compact.
This is the case for $V\subset\Lq$ and $W\subset\Cx0$.
We note that also the embeddings $W\subset V$, $V\subset H$
and $H\subset\Vp$ are compact.
Moreover, we often account for the \wk\ Poincar\'e inequality
\Beq
  \normaV v \leq C \Bigl( \normaH{\nabla v} + \bigl| \textstyle\iO v \bigr| \Bigr)
  \quad \hbox{for every $v\in V$}
  \label{poincare}
\Eeq
where $C$ depends only on~$\Omega$.
Furthermore, we repeatedly make use of the notation
\Beq
  Q_t := \Omega \times (0,t)
  \quad \hbox{for $t\in[0,T]$}
  \label{defQt}
\Eeq
and of \wk\ inequalities, namely, the \holder\ inequality and the elementary Young inequality:
\Beq
  ab \leq \delta a^2 + \frac 1{4\delta} \, b^2
  \quad \hbox{for every $a,b\geq 0$ and $\delta>0$}.
  \label{eleyoung}
\Eeq
Next, we introduce a tool that is generally used
in the context of problems related to the Cahn-Hilliard equations.
We define $\calN:\dom\calN\subset\Vp\to V$ as follows{\pier :}
\Beq
  \dom\calN := \graffe{\vstar\in\Vp: \ \<\vstar,1> = 0}
  \label{defdomN}
\Eeq
{\pier and,}
\Bsist
  &&\hbox{for $\vstar\in\dom\calN$,}\quad
  \hbox{$\calN\vstar$ is the unique $w\in V$ such that}
  \non
  \\
  && \iO \nabla w \cdot \nabla v = \< \vstar , v >
  \quad \hbox{for every $v\in V$}
  \aand
  \iO w = 0 \,.
  \label{defN}
\Esist
Note that problem \eqref{defN} actually has a unique solution $w\in V$
since $\Omega$ is also connected
and that $w$ solves the homogeneous Neumann problem for $-\Delta w=\vstar$
in the special case $\vstar\in H$.
It is easily checked that
\Bsist
  \hskip-1.5cm&& \< \zstar , \calN \vstar >
  = \< \vstar , \calN \zstar >
  = \iO (\nabla\calN\zstar) \cdot (\nabla\calN\vstar)
  \quad \hbox{for $\zstar,\vstar\in\dom\calN$}
  \label{simmN}
  \\
  \hskip-1.5cm&& \vstar \mapsto \normaVp\vstar
  := \normaH{\nabla\calN\vstar}
  = \< \vstar , \calN \vstar >^{1/2}
  \quad \hbox{{\pier is a norm on~$\dom\calN$}}
  \qquad
  \label{defnormaVp}
  \\
  \hskip-1.5cm&& \norma\vstar_{\Vp}
  \leq C \normaVp\vstar
  \quad \hbox{for some constant $C$ and every $\vstar\in\dom\calN$}
  \label{equivduale}
  \\
  \hskip-1.5cm&& 2 \< \dt\vstar(t) , \calN\vstar(t) >
  = \frac d{dt} \iO |\nabla\calN\vstar(t)|^2
  = \frac d{dt} \, \normaVp{\vstar(t)}^2
  \quad \aat \nonumber
  \\
  \hskip-1.5cm&&\quad \hbox{{\pier and for every $\vstar\in\H1\Vp$ satisfying
  $\<\vstar(t),1>=0$ for all $t\in[0,T]$.}}
  \label{dtcalN}
\Esist
Inequality~\eqref{equivduale} (where one could check that $C=1$ actually is suitable)
essentially says that $\normaVp\cpto$ and
the usual dual norm $\norma\cpto_{\Vp}$ are equivalent on~$\dom\calN$
(the opposite inequality is \sfw, indeed).
Finally, throughout the paper,
we use a small-case italic $c$ without any subscript for different constants,
that may only depend on~$\Omega$, the constant~$\gamma$,
the shape of the nonlinearities $p$, $\beta$ and~$\lambda$,
and the norms of the initial data related to assumption~\eqref{hpdati}.
We point out that $c$ does not depend on $\alpha$ nor on the final time~$T$
nor on the parameter $\eps$ we introduce in a forthcoming section.
For any parameter~$\delta$,
a~notation like~$c_\delta$ or $c(\delta)$ signals a constant that depends also on the parameter~$\delta$.
This holds, in particular, if $\delta$ is either $\a$, or $T$, or the pair $(\alpha,T)$.
The reader should keep in mind that the meaning of $c$ and $c_\delta$ might
change from line to line and even in the same chain of inequalities,
whereas those constants we need to refer to are always denoted by different symbols,
e.g., by a capital letter like in~\eqref{sobolev} or~by a letter with a proper subscript
as in~\eqref{classico}.


\section{Uniqueness}
\label{UniquenessSect}
\setcounter{equation}{0}

In this section, we prove the uniqueness part of Theorem~\ref{Wellposedness}, that is,
we pick two solutions $(\mu_i,u_i,\sigma_i,\xi_i)$, $i=1,2$, and show that they are the same.
As both $\a$ and $T$ are fixed, we avoid to stress the dependence of the constants on such parameters.
Moreover, as the solutions we are considering are fixed as well,
we can allow the values of $c$ to depend on them, in addition.
So, we write \eqref{unomenotre} and some of the equations of \Pbl\ for both solutions
and take the difference.
If we set $\mu:=\mu_1-\mu_2$ for brevity and analogously define $u$, $\sigma$ and~$\xi$,
we~have
\Bsist
  && \dt \bigl( \a \mu + u + \sigma \bigr)
  - \Delta (\mu + \sigma)
  = 0
  \label{diffunomenotre}
  \\
  && \mu = \a\dt u - \Delta u + \xi + \lambda'(u_1) - \lambda'(u_2)
  \label{diffseconda}
  \\
  && \dt\sigma - \Delta\sigma
  = R_2 - R_1 {\rev .}
  \label{diffterza}
\Esist
We note that \eqref{diffunomenotre} implies
$\iO(\a\mu+u+\sigma)(t)=\iO(\a\mu+u+\sigma)(0)=0$ for every~$t$,
so that $\calN(\a\mu+u+\sigma)$ is meaningful and we can use it to test~\eqref{diffunomenotre}. {\pier Then, we test \eqref{diffseconda}
by~$-u$ and sum the resulting equalities. Using}
the properties \eqref{defN}, \eqref{defnormaVp} and \eqref{dtcalN} of~$\calN$,
we have for every $t\in[0,T]$
\Bsist
  && \frac 12 \normaVp{\a\mu(t)+u(t)+\sigma(t)}^2
  + \intQt (\mu+\sigma) (\a\mu+u+\sigma)
  \non
  \\
  && \quad {}
  - \intQt \mu u
  + \frac \a 2 \iO |u(t)|^2
  + \intQt |\nabla u|^2
  + \intQt \xi u
  + \intQt \bigl( \lambda'(u_1) - \lambda'(u_2) \bigr) \, u
  = 0
  \non
\Esist
We note that the term involving $\xi$ is nonnegative since $\beta$ is monotone.
So, we rearrange and estimate the \rhs\ we obtain
by accounting for \eqref{equivduale} and the \Lip\ continuity of $\lambda'$
as follows
\Bsist
  && \frac 12 \normaVp{(\a\mu+u+\sigma)(t)}^2
  + \alpha \intQt |\mu|^2
  + \frac \a 2 \iO |u(t)|^2
  + \intQt |\nabla u|^2
  \non
  \\
  && \leq -\intQt \mu \sigma \, 
  {}{\rev -}{} \iot \< (\a\mu+u+\sigma)(s) , \sigma(s) > \, ds
  + c \intQt |u|^2
  \non
  \\
  && \leq  \delta \intQt |\mu|^2
  + \cd \intQt |\sigma|^2
  + c \intQt |u|^2
  \non
  \\
  && \quad {}
  + \delta \iot \normaV{\sigma(s)}^2 \, ds
  + \cd \iot \normaVp{(\a\mu+u+\sigma)(s)}^2 \, ds
  \label{stimadiffunomenotre}
\Esist
for every $\delta>0$.
Next, we test \eqref{diffterza} by $\sigma$ and get
\Bsist
  && \frac 12 \iO |\sigma(t)|^2
  + \intQt |\nabla\sigma|^2
  = \intQt (R_2 - R_1) \sigma
  \non
  \\
  && = \intQt \bigl( p(u_2) - p(u_1) \bigr) (\sigma_2 - \gamma\mu_2) \sigma
  - \intQt p(u_1) (\sigma - \gamma\mu) \sigma
  \non
  \\
  && \leq c \intQt |\sigma_2 - \gamma\mu_2| \, |u| \, |\sigma|
  + c \intQt |\sigma - \gamma\mu| \, |\sigma|
  \label{stimadiffterze}
\Esist
since $p$ is \Lip\ continuous and bounded.
Now, we estimate the last two integrals, separately.
By the regularity \eqref{regcomponenti} of $\mu_2$ and~$\sigma_2$
and the Sobolev inequality~\eqref{sobolev},
we~have
\Bsist
  && \intQt |\sigma_2 - \gamma\mu_2| \, |u| \, |\sigma|
  \leq \iot \norma{(\sigma_2 - \gamma\mu_2)(s)}_4 \, \norma{u(s)}_2 \, \norma{\sigma(s)}_4 \, ds
  \non
  \\
  && \leq c \iot \norma{u(s)}_2 \, \normaV{\sigma(s)} \, ds
  \leq \delta \intQt \bigl( |\sigma|^2 + |\nabla\sigma|^2 \bigr)
  + \cd \intQt |u|^2
  \non
\Esist
for every $\delta>0$.
On the other hand
\Beq
  \intQt |\sigma - \gamma\mu| \, |\sigma|
  \leq \intQt |\sigma|^2
  + c \intQt |\mu| \, |\sigma|
  \leq \delta \intQt |\mu|^2
  + \cd \intQt |\sigma|^2 .
  \non
\Eeq
At this point, we combine the last two estimates with \eqref{stimadiffterze}
and sum to~\eqref{stimadiffunomenotre}.
Then, we take $\delta$ small enough in order to absorb
the corresponding terms {\pier in} the \lhs.
By applying the Gronwall lemma,
we obtain $\mu=0$, $u=0$ and $\sigma=0$.
By comparison in \eqref{diffseconda} we also deduce $\xi=0$,
and the proof is complete.


\section{A priori estimates}
\label{Estimates}
\setcounter{equation}{0}

In this section, we introduce an approximating problem
and prove a number of a~priori estimates on its solution.
Some of the bounds we find may depend on $\a$ and~$T$,
while other ones are independent of such parameters.
The notation we use follows the general rule
explained at the end of Section~\ref{Statement}.
The estimates we obtain will be used in the {\danielle subsequent} sections
in order to prove our results.

For $\eps\in(0,1)$, the {\pier approximation to problem \accorpa{regcomponenti}{cauchy} is obtained} by simply replacing \eqref{regxi}~by
\Beq
  \xi = \betaeps(u)
  \label{defxieps}
\Eeq
where $\betaeps$ and the related functions $\Betaeps$ and $\Weps$
are defined on the whole of $\erre$ as follows
\Beq
  \Betaeps (r) :=
  \min_{s\in\erre} \Bigl( {\textstyle \frac 1 {2\eps}} \, (s-r)^2 + \Beta(s) \Bigr) , \quad
  \betaeps(r) := \frac d {dr} \Betaeps(r) , \quad
  \Weps(r) := \Betaeps(r) + {\pier \lambda}(r) .
  \label{defBbetaeps}
\Eeq
It turns out that $\Betaeps$ is a well-defined $C^1$ function
and that~$\betaeps$, the Yosida \regulariz ation of~$\beta$, is \Lip\ continuous.
Moreover the following properties
\Beq
  \hbox{$0\leq\Betaeps(r)\leq\Beta(r)$\aand $\Betaeps(r)\tosn\Beta(r)$ monotonically as $\eps\todx0$}
  \label{propYosida}
\Eeq
hold true for every $r\in\erre$ (see, e.g., \cite[Prop.~2.11, p.~39]{Brezis}).
Our starting point is the result stated below.

\Bprop
\label{Wellposednesseps}
Under the assumptions of Theorem~\ref{Wellposedness},
the approximating problem has a unique global solution.
\Eprop

Uniqueness {\pier is} already proved as a {\danielle special case} of the uniqueness part of Theorem~\ref{Wellposedness}.
As far as existence is concerned,
we avoid a detailed proof and just say that a Faedo-Galerkin method
(obtained by taking a base of~$V$, e.g., the base of the eigenfunctions of the Laplace operator
with homogeneous Neumann boundary conditions)
and {\pier some}~a~priori estimates very close
to the ones we are going to perform in the present section
{\pier lead} to the existence of a solution.
We also remark that the system of ordinary differential equations given by the Faedo-Galerkin scheme
has a unique global solution since $\betaeps$ is \Lip\ continuous
and the function $(\mu,u,\sigma)\mapsto p(u)(\sigma-\gamma\mu)$ on $\erre^3$
is smooth (since $p$ is~so) and sublinear (since $p$ is bounded).

{\pier From now on,} $(\mu,u,\sigma)=(\mueps,\ueps,\sigmaeps)$
stands for the solution to the approximating problem.
Accordingly, we define $\Reps$ by~\eqref{prima}.
However, we explicitely write the subscript $\eps$ only at the end of each calculation, for brevity.

\step
First a priori estimate

We multiply \eqref{prima} by $\mu$, {\rev\eqref{seconda} by $-\dt u$
and \eqref{terza} by $\sigma/\gamma$. Then, we {\pier add} all the equalities we obtain to each other,} integrate over~$Q_t$,
where $t\in(0,T)$ is arbitrary, and {\rev we} take advantage of the boundary conditions~\eqref{bc}.
We have
\Bsist
  && \frac \a 2 \iO |\mu(t)|^2
  - \frac \a 2 \iO |\muz|^2
  + \intQt \dt u \, \mu
  + \intQt |\nabla\mu|^2
  \non
  \\
  && {} + \frac 1 {2\gamma} \iO |\sigma(t)|^2
  - \frac 1 {2\gamma} \iO |\sigmaz|^2
  + \frac 1 \gamma \intQt |\nabla\sigma|^2
  {\pier {}- \intQt \mu \dt u
  + \a \intQt |\dt u|^2 }
  \non
  \\
  && {}
  + \frac 12 \iO |\nabla u(t)|^2
  - \frac 12 \iO |\nabla\uz|^2
  + \iO \Weps(u(t))
  - \iO \Weps(\uz)
  \non
  \\
  && = \intQt \bigl( R \mu - R (\sigma/\gamma) \bigr)
  = - \intQt p(u) \bigl( \gamma^{1/2} \mu - \gamma^{-1/2} \sigma \bigr)^2
  \label{pier1}
\Esist
and notice that two integrals cancel out. {\pier Moreover, we point out that
$$\iO \Weps(\uz) \leq \iO\Beta (\uz)+ \iO \lambda (\uz) \leq c $$
thanks to \accorpa{defBbetaeps}{propYosida}, \eqref{hpdati} and \eqref{hplambda}.}
By rearranging {\pier in \eqref{pier1}} and using \eqref{hpdati} and~\eqref{propYosida}, we thus deduce
\Bsist
  \hskip-.8cm&& \frac \a 2 \iO |\mu(t)|^2
  + \intQt |\nabla\mu|^2
  + \frac 1 {2\gamma} \iO |\sigma(t)|^2
  + \frac 1 \gamma \intQt |\nabla\sigma|^2
  \non
  \\
  \hskip-.8cm&& {} + \a \intQt |\dt u|^2
  + \frac 12 \iO |\nabla u(t)|^2
  + \iO \Weps(u(t))
  + \intQt p(u) \bigl( \gamma^{1/2} \mu - \gamma^{-1/2} \sigma \bigr)^2
  \leq c \,.
  \qquad
  \label{perprimastima}
\Esist
On the other hand, by simply integrating \eqref{unomenotre} over $\Omega$ \aat\ and using~\eqref{bc},
we~get
\Beq
  \iO \bigl( \a\mu(t) + u(t) + \sigma(t) \Bigr)
  = \iO \bigl( \a\muz + \uz + \sigmaz \Bigr)
  = c
  \non
\Eeq
whence immediately (here $|\Omega|$ is the Lebesgue measure of~$\Omega$)
\Bsist
  && \Bigl| \iO u(t) \Bigr|^2
  \leq \bigl( c + \a \norma{\mu(t)}_1 + \norma{\sigma(t)}_1 \bigr)^2
  \leq 3c^2 + 3\a^2 \norma{\mu(t)}_1^2 + 3 \norma{\sigma(t)}_1^2
  \non
  \\
  && \leq c + 3|\Omega| \bigl( \a^2 \norma{\mu(t)}_2^2 + \norma{\sigma(t)}_2^2 \bigr)
  \leq c + D \bigl( \a \norma{\mu(t)}_2^2 + {\textstyle\frac 1\gamma} \norma{\sigma(t)}_2^2 \bigr)
  \non
\Esist
where $D:=3|\Omega| \max\{1,\gamma\}$.
At this point, we multiply the above inequality by $1/(4D)$,
sum what we obtain to~\eqref{perprimastima} and rearrange.
{\pier Using also} the Poincar\'e inequality~\eqref{poincare},
we conclude~that
\Bsist
  && \a^{1/2}\norma\mueps_{\L\infty H}
  + \norma{\nabla\mueps}_{\L2H}
  \non
  \\
  && {} \quad
  + \a^{1/2} \norma{\dt\ueps}_{\L2H}
  + \norma\ueps_{\L\infty V}
  + \norma{\Weps(\ueps)}_{\L\infty\Luno}
  \non
  \\
  && {} \quad
  + \norma\sigmaeps_{\L\infty H}
  + \norma{\nabla\sigmaeps}_{\L2H}
  \non
  \\
  && {} \quad
  + \norma{(p(\ueps))^{1/2}(\gamma^{1/2}\mu-\gamma^{-1/2}\sigma)}_{\L2H}
  \leq c \,.
  \qquad
  \label{primastima}
\Esist
By recalling the definition {\pier of~$\Reps$ in~\eqref{prima}}
(this is the notation for the approximating problem, indeed),
we~have
\Beq
  \Reps = p(\ueps) (\sigmaeps - \gamma\mueps)
  = -(p(\ueps))^{1/2} \gamma^{1/2} \cdot (p(\ueps))^{1/2}(\gamma^{1/2}\mueps-\gamma^{-1/2}\sigmaeps).
  \non
\Eeq
As $p$ is bounded, {\pier from \eqref{primastima} we infer that}
\Beq
  \norma\Reps_{\L2H} \leq c \,.
  \label{stimaR}
\Eeq

\step
Second a priori estimate

We write \eqref{prima} in the form
$\a\dt\mu-\Delta\mu=R-\dt u$
and multiply this equation by $\dt\mu$.
By integrating over~$Q_t$, we obtain
\Bsist
  && \a \intQt |\dt\mu|^2
  + \frac 12 \iO |\nabla\mu(t)|^2
  - \frac 12 \iO |\nabla\muz|^2
  \non
  \\
  && = \intQt (R - \dt u) \dt\mu
  \leq \norma{R-\dt u}_{\L2H} \norma{\dt\mu}_{\L2H} \,.
  \non
\Esist
{\pier Thanks} to \eqref{primastima} and \eqref{stimaR},
we conclude that {\pier
\Beq
  \norma{\dt\mueps}_{\L2H}
  + \norma\ueps_{\L\infty V}
  \leq \ca \,.
  \label{secondastima}
\Eeq
}%
In the same way, we derive the analogue for $\sigmaeps$ by using~\eqref{terza}:
{\pier however, in this case we obtain the better estimate
\Beq
  \intQt |\dt\sigma|^2
  + \frac 12 \iO |\nabla\sigma(t)|^2
  \leq  \frac 12 \iO |\nabla\sigmaz|^2
  + \norma{R}_{\L2H} \norma{\dt\sigma}_{\L2H} \,.
  \non
\Eeq
which implies, along with \eqref{primastima},
\Beq
\norma{\dt\sigmaeps}_{\L2H}  + \norma\sigmaeps_{\L\infty V} \leq c .
\label{pier2}
\Eeq
}%

\step
Third a priori estimate

By writing our system in the form
\Beq
  \a \dt\mu - \Delta\mu = R - \dt u , \quad
  -\Delta u + \betaeps(u) = \mu - \a\dt u - \lambda'(u), \quad
  \dt\sigma - \Delta\sigma = -R
  \non
\Eeq
and accounting for the bounds already found,
it is \sfw\ to~derive the following estimates
(multiply by $-\Delta v$ with $v=\mu,u,\sigma$, respectively,
and {\rev use $\betaeps'\geq0$})
\Beq
  \norma{\Delta\mu}_{\L2H} \leq \ca \, \quad
  \norma{\Delta u}_{\L2H} \leq \caT \, \quad
  {\pier  \norma{\Delta\sigma}_{\L2H} \leq c} \,.
  \label{stimeDelta}
\Eeq
By the elliptic regularity theory, {\rev \eqref{seconda} and the boundary conditions \eqref{bc}},
we conclude~that
\Bsist
  && \norma\mueps_{\L2W} \leq \caT
  \aand
  \norma\sigmaeps_{\L2W} \leq \cT
  \label{stimaWmusigma}
  \\
  && \norma\ueps_{\L2W} + \norma{\betaeps(\ueps)}_{\L2H} \leq \caT \,.
  \label{stimaWu}
\Esist


\section{Existence on a finite time interval}
\label{Existence}
\setcounter{equation}{0}

In this section, we conclude the proof of Theorem~\ref{Wellposedness}
by showing the existence of a solution.
As $\a$ and $T$ are fixed now, we can account for all of the estimates of the previous section.
{\pier Owing to standard weak compactness {\pier arguments} as well as by the strong compactness result in ~\cite[Sect.~8, Cor.~4]{Simon}, it turns out that}
the following convergence holds
\Bsist
  && \mueps \to \mu , \quad
  \ueps \to u , \quad
  \sigmaeps \to \sigma
  \non
  \\
  && \qquad \hbox{weakly {\pier star} in $\H1H{\pier{}\cap\L{\infty}V{}}\cap\L2W$,}
  \non
  \\
  && \qquad \hbox{strongly in ${\pier C^0([0,T];H)\cap{}}\L2V$ and \aeQ}
  \label{convsoluzeps}
  \\
  && \Reps \to \rho
  \aand
  \betaeps(\ueps) \to \xi
  \quad \hbox{weakly in $\L2H$}
  \label{convRxieps}
\Esist
at least for a subsequence.
{\pier The strong} convergence in $\C0H$ {\pier entails that}
the limiting functions satisfy the initial conditions~\eqref{cauchy}.
Moreover, the {\pier pointwise} convergence~a.e.\ {\rev and assumption 
\eqref{hpp} imply} that
$\Reps$ converges to $R:=p(u)(\sigma-\gamma\mu)$ \aeQ,
so that $\rho=R$.
Furthermore, as $\lambda'$ is \Lip\ continuous,
$\lambda'(\ueps)$ converges to $\lambda'(u)$ strongly in $\L2H$.
Finally, a~standard monotonicity argument
(see, e.g., \cite[Lemma~1.3, p.~42]{Barbu})
{\pier based on the weak convergence $\ueps \to u$,
$\betaeps(\ueps) \to \xi$ in $L^2 (Q)$ and on the property
{\gianni
\Beq
  \limsup_{\eps\searrow 0} \int_Q \betaeps(\ueps) \ueps
  {\rev {}=\lim_{\eps\searrow 0} \int_Q \betaeps(\ueps) \ueps  
  ={}} \int_Q \xi u
  \non
\Eeq
}}%
{\rev (easily following from \accorpa{convsoluzeps}{convRxieps})} 
yields $\xi\in\beta(u)$ \aeQ.
Therefore{\pier , the quadruplet} $(\mu,u,\sigma,\xi)$
satisfies \Regsoluz\ and {\pier solves \Pbl}.

We conclude this section by {\pier recovering the uniform estimates
for the solution $(\mu,u,\sigma,\xi)$ to the problem \Tuttopbl}.
First of all, we can speak of a unique solution
on the time half line~$[0,+\infty)$.
For such a solution, {\pier the estimates we found
for the approximating problem still hold},
i.e., we~have
\Bsist
  && \a^{1/2}\norma\mu_{\LL\infty H}
  + \norma{\nabla\mu}_{\LL2H}
  \non
  \\
  && {} \quad
  + \a^{1/2} \norma{\dt u}_{\LL2H}
  + \norma u_{\LL\infty V}
  + \norma{\calW(u)}_{\LL\infty\Luno}
  \non
  \\
  && {} \quad
  {\pier{}+ \norma{\dt\sigma}_{\LL2H}}
  + \norma\sigma_{{\pier \LL\infty V}}
  + \norma{\nabla\sigma}_{\LL2H}
  + \norma{\Delta\sigma}_{\LL2H}
  \qquad
  \non
  \\
  && {} \quad
  + \norma R_{\LL2H}
  \leq c
  \label{daprimastima}
  \\
  && {\pier \norma{\dt\mu}_{\LL2H}
  +{}} \norma\mu_{\LL\infty V}
  + \norma u_{\LL\infty V}
  + \norma\sigma_{\LL\infty V}
  \leq \ca
  \label{dasecondastima}
  \\
  && \norma\mu_{\L2W}
  + \norma u_{\L2W}
  + \norma\sigma_{\L2W}
  + \norma\xi_{\L2H}
  \leq \caT
  \label{dastimeW}
\Esist
the last one for every $T\in(0,+\infty)$.
This simply follows from the weak semicontinuity of the norms
for all the inequalities but the one involving $\calW(u)$.
As far as the latter is concerned,
{\pier we note that, for all $t\geq 0$, $\ueps (t)\to u(t)$ strongly in~$H$.
Then, using \eqref{hplambda} and the
mean value theorem, it is not difficult to check that
\Beq
  |\lambda(\ueps(t)) - \lambda(u(t))|
  \leq c \, |\ueps(t) - u(t)| \left( 1 + |\ueps(t)| + |u(t)| \right)
  \non
\Eeq
and consequently $\lambda(\ueps(t))$ converges to $\lambda(u(t))$ strongly in $\Lx1$.
Hence, in view of \eqref{defWbeta} and \eqref{defBbetaeps}
it suffices to prove that
\Beq
  \iO \Beta(u(t))
  \leq \liminf_{\eps\todx 0} \iO \Betaeps(\ueps(t)) .
  \label{pier3}
\Eeq
To this end, we fix $\eps'>0$ for a while.
By accounting for the lower semicontinuity of $\Beta_{\eps'}$ and the inequality
$\Beta_{\eps'}(s)\leq\Betaeps(s)$ which holds for every $s\in\erre$ and
$\eps\in(0,\eps')$ (see~\eqref{defBbetaeps}), we~obtain
\Beq
  \iO \Beta_{\eps'}(u(t))
  \leq \liminf_{\eps\todx 0} \iO \Beta_{\eps'}(\ueps(t))
  \leq \liminf_{\eps\todx 0} \iO \Betaeps(\ueps(t)) .
  \label{pier4}
\Eeq
Now, we let $\eps'$ vary and recall \eqref{propYosida} in terms of $\eps'$.
Thus, the Beppo Levi monotone convergence theorem implies that
\Beq
  \iO \Beta(u(t))
  = \lim_{\eps'\todx 0} \iO \Beta_{\eps'}(u(t))
  \label{pier5}
\Eeq
and combining \eqref{pier5} and \eqref{pier4} yields~\eqref{pier3}.
Therefore, \eqref{daprimastima} follows.}


\section{Long time \bhv}
\label{LongtimeSect}
\setcounter{equation}{0}

In this section, we prove Theorem~\ref{Longtime}.
From \eqref{dasecondastima} we see that the {\rev omega limit} $\omega$ we are interested in is non-empty.
It remains to \characteriz e its elements as in the statement.
So, we fix $(\muo,\uo,\sigmao)\in\omega$ and a sequence $\{\tn\}$
according to definition~\eqref{omegalimit}.
We set for convenience
\Beq
  \vn(t) := v(t+\tn)
  \quad \hbox{for $t\geq0$}
  \quad \hbox{with} \quad
  v = \mu ,\, u ,\, \sigma ,\, \xi ,\, R
  \label{defcode}
\Eeq
and study the \bhv\ of such functions in a fixed finite time interval~$(0,T)$.
First of all, we notice that the quadruplet $(\mun,\un,\sigman,\xin)$
solves the problem obtained from problem \Pbl\
by replacing the initial condition~\eqref{cauchy}
by the following one
\Beq
  \mun(0) = \mu(\tn) , \quad
  \un(0) = u(\tn)
  \aand
  \sigman(0) = \sigma(\tn) .
  \label{cauchyn}
\Eeq
On the other hand, by~\eqref{dasecondastima}, the new initial data are bounded in~$V$
and {\pier\eqref{daprimastima}} provides a uniform $L^1$-estimate for $\calW(\un(0))=\calW(u(\tn))$.
Therefore, the dependence of the constants
on the norms of the initial data just mentioned
leads to a dependence only on the norms involved in our assumptions~\eqref{hpdati}.
Moreover, we observe that, for every Banach space~$Z$,
if some function $v$ belongs to $\LL2Z$
and $\vn$ is related to $v$ as in~\eqref{defcode},
we trivially have
\Beq
  \lim_{n\to\infty} {\pier {}\int_0^T} \norma{\vn(t)}_Z^2 \, dt
  {\pier {}\leq{}} \lim_{n\to\infty} \int_{\tn}^\infty \norma{v(t)}_Z^2 \, dt
  = 0
  \non
\Eeq
so that $\vn\to0$ strongly in $\L2Z$.
At this point, we can derive the estimates and the convergence we need
from \accorpa{daprimastima}{dastimeW}.
We {\pier infer that}
\Bsist
  && \norma\mun_{\LL\infty V}
  + \norma\un_{\LL\infty V}
  + \norma\sigman_{\LL\infty V}
  \leq \ca
  \non
  \\
  && \norma\mun_{\L2W}
  + \norma\un_{\L2W}
  + \norma\sigman_{\L2W}
  + \norma\xin_{\L2H}
  \leq \caT
  \non
  \\
  && \nabla\mun \to 0
  \aand
  \nabla\sigman \to 0
  \quad \hbox{strongly in $(\L2H)^3$}
  \non
  \\
  && \dt\mun \to 0 , \quad
  \dt\un \to 0 , \quad
  \dt\sigman \to 0
  \aand
  \Rn \to 0
  \quad \hbox{strongly in $\L2H$} .
  \non
\Esist
Therefore, at least for a subsequence, we also have
\Bsist
  && \mun \to \mui , \quad
  \un \to \ui
  \aand
  \sigman \to \sigmai
  \non
  \\
  && \qquad  \hbox{weakly star in $\L\infty V\cap \L2W$}
  \label{convcode}
  \\
  && \mun \to \mui \, \quad
  \un \to \ui
  \aand
  \sigman \to \sigmai
  \quad \hbox{strongly in $\L2H$}
  \label{strongcode}
  \\
  && \xin \to \xii
  \quad \hbox{weakly in $\L2H$}
  \label{convxin}
\Esist
the strong convergence \eqref{strongcode}
being a consequence of \eqref{convcode} and of the bounds for the time derivatives.
In particular, {\pier thanks to the \Lip\ continuity of~$p$,}
we derive that $\Rn$ converges to $p(\ui)(\sigmai-\gamma\mui)$
{\pier strongly} in~$L^1(Q)$, whence
\Beq
  p(\ui)(\sigmai-\gamma\mui) = 0
  \quad \aeQ .
  \label{relazi}
\Eeq
Furthermore, $\mui$ and $\sigmai$ are constant functions and $\ui$~is time independent.
We denote the constant values of $\mui$ and $\sigmai$ by $\mus$ and~$\sigmas$, respectively,
and set $\us:=\ui(t)$ for $t\in(0,T)$.
By taking the limit in~\eqref{seconda} written for $\mun$ and~$\un$,
we see that the pair $(\ui,\xii)$ solves the following problem
\Beq
  \mus = -\Delta\ui + \xii + \lambda'(\ui)
  \quad \aeQ
  \aand
  {\pier \dn\ui = 0
  \quad \aeS}.
  \non
\Eeq
In particular, $\xii$~also is time independent, $\xii(t)=\xis$ for $t\in(0,T)$,
and the above boundary value problem and \eqref{relazi} become
\Bsist
  \mus = -\Delta\us + \xis + \lambda'(\us)
  \aand
  p(\us)(\sigmas-\gamma\mus) = 0
  \quad \aeO , \quad\non
  \\
  {\pier \dn\us = 0 \quad \aeG}.
  \label{relazs}
\Esist
Now, as in Section~\ref{Existence},
we derive both $\xii\in\beta(\ui) \ {\pier \aeQ} $, i.e., $\xis\in\beta(\us)\ {\pier\aeO}$,  and the convergence
\Beq
  \mun \to \mui \,, \quad
  \un \to \ui, \quad
  {\pier {}\sigman \to \sigmai
  \quad \hbox{strongly in $\C0H\cap\L2V$}}.
  \non
\Eeq
It follows that $\mu(\tn)=\mun(0)$ converges to $\mui(0)=\mus$ {\pier in~$H$.
As $\mu(\tn)$} converges to $\muo$ in~$V$,
we infer that $\muo=\mus$.
In the same way we obtain $\uo=\us$ and $\sigmao=\sigmas$.
Therefore, we also have from~\eqref{relazs}
\Beq
  \muo \in -\Delta\uo + \beta(\uo) + \lambda'(\uo)
  \aand
  p(\uo)(\sigmas-\gamma\mus) = 0
  \quad \aeO
  \non
\Eeq
and the proof is complete.

\Brem
\label{Altrospazio}
Even though we have to confine ourselves to study the {\rev omega limit}
of an initial datum satisfying~\eqref{hpdati},
we could take a phase space $\Phi$ that is larger than \eqref{phasespace}
and is endowed with a weaker topology.
This may lead to further properties of~$\omega$.
For instance, if we choose
$\Phi=(\Ldue)^3$ with the strong topology,
estimate~\eqref{daprimastima} implies that
the whole trajectory of the initial datum is relatively compact in~$\Phi$,
so that general results (see, e.g., \cite[Lemma 6.3.2, p.~239]{Zheng})
ensure that $\omega$ is invariant, compact and connected in the $L^2$ topology.
\Erem


\section{Asymptotics and limit problem}
\label{Aphazero}
\setcounter{equation}{0}

In this section, we perform the proof of Theorem~\ref{Asymptotics}.
As $T$ is fixed, we avoid stressing the dependence of the constants on~$T$.

$i)$~As in the statement,
$(\mua,\ua,\sigmaa,\xia)$ (where $\xia={\pier \beta} (\ua)$ since $\pier \Beta$ is smooth)
is~the solution to problem \Pbl\ and we define $\Ra$ accordingly.
We recall that \eqref{daprimastima} implies
$\norma\ua_{\L\infty V}\leq c$,
whence also $\norma\ua_{\L\infty{\Lx6}}\leq c$
due to the Sobolev inequality~\eqref{sobolev}.
By {\pier the} assumption~\eqref{classico}, we infer that
\Beq
  \norma{\beta(\ua)+{} {\pier\lambda}'(\ua)}_{\L\infty H}
  = \norma{\calW'(\ua)}_{\L\infty H}
  \leq c \bigl( \norma\ua_{\L\infty{\Lx6}}^3 + 1 \bigr)
  \leq c \,.
  \label{stimaWp}
\Eeq
Now, we integrate \eqref{seconda} over~$\Omega$
and use the homogeneous Neumann boundary condition for~$\rev \ua$.
Then, we square and integrate over $(0,T)$ with respect to time.
We obtain
\Bsist
  && \ioT \Bigl| \iO \mua(t) \Bigr|^2 \, dt
  = \ioT \Bigl| \iO \bigl( \a\, {\pier \dt}{}\ua(t) + \calW'(\ua(t)) \bigr) \Bigr|^2 \, dt
  \non
  \\
  && \leq 2 \,{\pier |\Omega| }\,  \a^2 \intQ |\dt\ua|^2
  + 2  \, {\pier |\Omega| } \intQ |\calW'(\ua)|^2
  \leq c
  \non
\Esist
the last inequality {\pier following from} \eqref{daprimastima} and~\eqref{stimaWp}.
{\pier Then, recalling the estimate for $\nabla\mua$ in~\eqref{daprimastima}}
and owing to the Poincar\'e inequality~\eqref{poincare},
we conclude~that
\Beq
  \norma\mua_{\L2V} \leq c \,.
  \label{stimamu}
\Eeq
By comparison in~\eqref{seconda}, we deduce $\norma{\Delta\ua}_{\pier \L2H}\leq c$,
whence also
\Beq
  \norma\ua_{\L2W} \leq c
  \label{stiamu}
\Eeq
by elliptic regularity.
Now, we test \eqref{prima} by an arbitrary $v\in\L2V$ and~get
\Beq
  \Bigl| \intQ \dt(\a\mua+\ua) \, v \Bigr|
  = \Bigl| \intQ \bigl( \Ra v - \nabla\mua \cdot \nabla v \bigr) \Bigr|
  \leq c \norma v_{\L2V}
  \non
\Eeq
by the estimate \eqref{daprimastima} of $\Ra$ and~\eqref{stimamu}.
This means~that
\Beq
  \norma{\dt(\a\mua+\ua)}_{\L2\Vp} \leq c \,.
  \label{stimadtVp}
\Eeq
Next, {\pier still from \eqref{daprimastima} it follows that}
\Beq
  \norma\sigmaa_{\H1H\cap\L2W} \leq c \,.
  \label{stimasigma}
\Eeq
At this point, we can use weak and weak star compactness and conclude that
\Bsist
  & \mua \to \mu
  & \hbox{weakly in $\pier \L2V$}
  \label{convmua}
  \\
  & \ua \to u
  & \hbox{weakly star in $\L\infty V\cap\L2W$}
  \label{convua}
  \\
  & \sigmaa \to \sigma
  & \hbox{weakly in $\H1H\cap\L2W$}
  \label{convsigmaa}
  \\
  & \dt(\a\mua + \ua) \to \zeta
  & \hbox{weakly in $\L2\Vp$}
  \label{convdt}
\Esist
at least for a subsequence.
This proves the part~$i)$ of the statement but~\eqref{convdtamuu},
which is more precise than~\eqref{convdt} and is justified in the next step.

$ii)$~Take any triplet $(\mu,u,\sigma)$ satisfying the above convergence
(note that \eqref{convdt} is weaker than \eqref{convdtamuu}):
we prove that it solves the limit problem \accorpa{primalim}{cauchylim}
and that $\zeta=\dt u$ (so~that \eqref{convdtamuu} holds).
First of all, we notice that \eqref{convsigmaa} implies $\sigma(0)=\sigmaz$.
Next, as \accorpa{convmua}{convua} imply that $\a\mua+\ua\to u$ weakly in~$\L2V$
and \eqref{convdt} holds,
we can apply the Lions-Aubin {\pier theorem}
(see, e.g., \cite[Thm.~5.1, p.~58]{Lions})
and deduce that
\Beq
  \a\mua + \ua \to u
  \quad \hbox{strongly in $\L2H$}.
  \non
\Eeq
On the other hand $\a\mua\to0$ in $\L2V$ by~\eqref{convmua}.
Hence
\Beq
  \ua \to u
  \quad \hbox{strongly in $\L2H$}
  \aand
  \zeta = \dt u .
  \label{forteua}
\Eeq
By {\pier a} standard argument (the same as in Section~\ref{Existence}),
we can {\pier identify} the limit of $\Ra$ as $p(u)(\sigma-\gamma\mu)$
and the limits of the other nonlinear terms{\pier . Thus, we}
conclude that $(\mu,u,\sigma)$ solves \accorpa{primalim}{bclim}
({\pier in~fact, one proves an equivalent} integrated version of \eqref{primalim}
rather than \eqref{primalim} itself).
It remains to check the first {\pier condition in}~\eqref{cauchylim}.
By also accounting for~\eqref{convdt}, we see that
$\a\mua+\ua$ converges to $u$ weakly in $\C0\Vp$.
This implies that
\Beq
  \a\muz + \uz
  = (\a\mua+\ua)(0) \to u(0)
  \quad \hbox{weakly in $\Vp$} .
  \non
\Eeq
On the other hand, $\a\muz+\uz\to\uz$ strongly in $V$.
Therefore, $u(0)=\uz$, and the proof of $ii)$ is complete.

$iii)$~A formal estimate that leads to \accorpa{regmubis}{regubis}
could be obtained by testing \eqref{primalim} by~$\dt u$,
differentiating \eqref{secondalim} with respect to time,
testing the obtained equality by~$\mu$ and {\rev adding} up.
{\rev Here we perform} the correct procedure, namely,
the discrete version of the formal one,
by introducing a time step $h\in(0,1)$.
For simplicity, we allow the (variable) value of the constant $c$
to depend on the {\gianni norm $\normaW\uz$} involved in~\eqref{hpdatibis}
and on the solution we are considering (which is fixed).
Of course, $c$~does not depend on~$h$.
{\gianni
First of all, we introduce a notation}.
For $v\in\Lm2H$ and $h\in(0,1)$, we define
the mean $\vh\in\L2H$ and the difference quotient $\dh v\in\L2H$
by setting for $t\in(0,T)$
\Beq
  \vh(t) := \frac 1h \int_{t{-}h}^t \!\!\! v(s) \, ds
  = \int_0^1 \! v(t{-}h\tau) \, d\tau
  \aand
  \dh v(t) := \dt \vh(t)
  = \frac {v(t)-v(t{-}h)} h
  \label{defvh}
\Eeq
and {\pier we do the same} if $v\in(\Lm2H)^3$ in order to treat gradients.
We notice that
\Beq
  \norma\vh_{\L2H} \leq \norma v_{\Lm2H}
  \label{stimavh}
\Eeq
as we show at once.
We have indeed
\Bsist
  && \norma\vh_{\L2H}^2
  {\pier {}\leq{}} \iO \ioT \int_0^1 |v(x,t-h\tau)|^2 \,d\tau \,dt \,dx
  = \int_0^1 \iO \ioT |v(x,t-h\tau)|^2 \,dt \,dx \,d\tau
  \non
  \\
  && = \int_0^1 \iO \int_{-h\tau}^{T-h\tau} |v(x,s)|^2 \, ds \,dx \,d\tau
  \leq \int_0^1 \iO \int_{-1}^T |v(x,s)|^2 \, ds \,dx \,d\tau
  = \norma v_{\Lm2H}^2 \,.
  \non
\Esist
{\gianni
As we are going to apply \eqref{stimavh} to $\mu$ and~$R$,
we need to extend such functions to the whole of $\Omega\times(-1,T)$.
Our tricky construction is based on assumption \eqref{hpdatibis} on~$\uz$
and involves $u$ as well.
We first solve a backward variational problem
with the help of the theory of linear abstract equations.
We~set
\Beq
  \Hz := \Bigl\{ v \in H : \ \textstyle\iO v =0 \Bigr\}
  \aand
  \Wz := W \cap \Hz
  \non
\Eeq
and construct the Hilbert triplet $(\Wz,\Hz,\Wzp)$,
where $\Wzp$ is the dual space of~$\Wz$,
by embedding $\Hz$ into $\Wzp$ in the standard way.
In the sequel, the symbol $\<\cpto,\cpto>$ denotes the duality pairing
between $\Wzp$ and~$\Wz$.
We introduce the continuous bilinear form $a:\Wz\times\Wz\to\erre$ by setting
\Beq
  a(z,v) := \iO (\Delta z) \, (\Delta v)
  \quad \hbox{for $z,v\in\Wz$}
  \non
\Eeq
and observe that $a(v,v)+\normaH v^2\geq\alpha\normaW v^2$
for some $\alpha>0$ and every $v\in\Wz$,
thanks to the elliptic regularity theory.
We also notice that $\Delta\uz\in\Hz$ since $\uz\in W$.
Therefore, as is well known
(e.g., \cite[Prop.~2.3 p.~112]{Showalter}),
there exists a unique $z$ satisfying
\Bsist
  & z \in \Hmz 1\Wzp \cap \Cmz 0\Hz \cap \Lmz2\Wz &
  \label{added1}
  \\
  & - \< \dt z(t) , v > + a(z(t),v) = 0
  \quad \hbox{for every $v\in\Wz$ and for a.a.\ $t\in(-1,0)$} & \qquad
  \label{added2}
  \\
  & z(0) = - \Delta\uz &
  \label{added3}
\Esist
and we also have
\Beq
  \norma z_{\Hmz 1\Wzp\cap\Lmz\infty H\cap\Lmz 2W}
  \leq c \normaH{z(0)}
  \leq c \normaW\uz
  = c \,.
  \label{regz}
\Eeq
As $z\in\Cmz 0\Hz$,
for every $t\in[-1,0]$ we have that $z(t)\in\dom\calN$ (see~\eqref{defdomN}).
Hence, we can define a function $w$ by setting
\Beq
  w(t) := \calN(z(t))
  \quad \hbox{for every $t\in[-1,0]$}
  \label{defw}
\Eeq
and it turns out that $w\in\Cmz 0W$:
the restriction of $\calN$ to $\Hz$
is an isomorphism from $\Hz$ onto~$\Wz$, indeed.
Moreover, $w$~is even smoother.
Namely, from~\eqref{regz} we have~that
\Beq
  \norma w_{\Lmz\infty W}
  \leq c
  \aand
  \norma{\dt w}_{\Lmz 2H}
  \leq c \,.
  \label{regw}
\Eeq
Here, the former is due to the above argument
and we now prove the latter.
Clearly, an estimate on the difference quotients is sufficient to conclude.
We observe that the operator $-\Delta:\Wz\to\Hz$
is a well-defined isomorphism.
Thus, the same property is enjoyed by its adjoint operator
$(-\Delta)^*:\Hz\to\Wzp$ given~by
\Beq
  \< (-\Delta)^* y , v > = \iO y (-\Delta v)
  \quad \hbox{for every $y\in\Hz$ and $v\in\Wz$} .
  \non
\Eeq
Hence, for every $\wstar\in\Wzp$ there exists a unique $y\in\Hz$
such that $(-\Delta)^*y=\wstar$,~i.e.,
\Beq
  \iO y (-\Delta v) = \< \wstar , v >
  \quad \hbox{for every $v\in\Wz$}
  \label{invaggiunto}
\Eeq
and the estimate
$\normaH y\leq C\norma\wstar_{\Wzp}$ holds true with $C$ depending only on~$\Omega$.
Assume now $h\in(0,1)$ and $t\in(-1{+}h,0)$.
From the definition \eqref{defw} of~$w$
we immediately derive that
$y=\dh w(t)$ belongs to $\Hz$ and satisfies
\Beq
  \iO y (-\Delta v)
  = \iO (-\Delta\dh w(t)) v
  = \iO (\dh z(t)) v
  \quad \hbox{for every $v\in\Wz$}
  \non
\Eeq
i.e., it fulfils \eqref{invaggiunto} with $\wstar=\dh z(t)$.
Therefore, we have
$\normaH{\dh w(t)}\leq c\norma{\dh z(t)}_{\Wzp}$
for every $t\in(-1{+}h,0)$,
and \eqref{regz} immediately implies
\Beq
  \int_{-1+h}^0 \normaH{\dh w(t)}^2 \, dt
  \leq c \int_{-1+h}^0 \norma{\dh z(t)}_{\Wzp}^2 \, dt
  \leq c
  \quad \hbox{for every $h\in(0,1)$}.
  \non
\Eeq
This proves the second estimate in~\eqref{regw}.
Once \eqref{regw} is completely established, we go on.
We term $(\uz)_\Omega$ the mean value of~$\uz$
and notice that $\uz-(\uz)_\Omega$ {\rev coincides with $w(0)$ 
given by \eqref{defw} and \eqref{added3}.
Therefore, 
we are suggested to define $u$ in $(-1,0)$ by setting}
\Beq
  u(t) := w(t) + (\uz)_\Omega
  \quad \hbox{for $t\in(-1,0)$}.
  \label{defprolu}
\Eeq
By doing that, we have both the estimates
\Beq
  \norma u_{\Lmz\infty W}
  \leq c
  \aand
  \norma{\dt u}_{\Lmz 2H}
  \leq c
  \label{stimeprolu}
\Eeq
and the fact that $u(t)\to\uz$ (e.g., in~$H$) as $t\tosn 0$.
This implies that the extended function $u\in\Lm2H$
(which is continuous in $[-1,0]$ and $[0,T]$, separately)
does not jump at $t=0$.
Thus, $u\in C^0([-1,T];H)$ and its time derivative $\dt u$
is a $\Vp$-valued function (rather than a distribution)
on~the whole of~$(-1,T)$,
so that \eqref{primalim} will hold in the whole of $(-1,T)$
whenever we properly extend $\mu$ and~$R$.
As the former is concerned, we~set
\Beq
  \mu(t) := z(t) + \calW'(u(t)) = -\Delta u(t) + \calW'(u(t))
  \quad \hbox{for $t\in(-1,0)$}.
  \label{defprolmu}
\Eeq
Now, in order to properly extend~$R$,
we check that $\Delta\mu$ is a well defined function.
Indeed, $\Delta z\in\Lmz 2H$ by \eqref{regz}.
On the other hand
\Beq
  \Delta\calW'(u)
  = \calW'''(u) |\nabla u|^2 + \calW''(u) \Delta u
  \in \Lmz\infty H
  \non
\Eeq
on account of \eqref{stimeprolu}, assumption \eqref{classico} on $\calW$
and the Sobolev embedding $W\subset\Wx{1,4}$.
Hence, we can~set
\Beq
  R(t) := \dt w(t) - \Delta\mu(t) = \dt u(t) - \Delta\mu(t)
  \quad \hbox{for $t\in(-1,0)$} .
  \label{defprolR}
\Eeq
We obtain
\Beq
  \norma\mu_{\Lmz 2W} \leq c
  \aand
  \norma R_{\Lmz 2H} \leq c \,.
  \label{stimeprolmuR}
\Eeq
Notice that it is confirmed that
\eqref{primalim} holds for {\gianni every $v\in V$ and a.e.\ in $(-1,T)$}.
Moreover, \eqref{secondalim} is satisfied a.e.\ in~$\Omega\times(-1,T)$
(while \eqref{defR} still holds only in~$Q$).
By collecting \eqref{stimeprolu}, \eqref{stimeprolmuR}
and the regularity of the original functions, we deduce, in particular, the following estimates
\Beq
  \norma u_{\Lm\infty V} \leq c \,, \quad
  \norma\mu_{\Lm2V} \leq c
  \aand
  \norma R_{\Lm2H} \leq c \,.
  \label{stimeprol}
\Eeq
At this point, we come back to \eqref{stimeprol}
and apply it to the extended $\mu$ and~$R$.
By also observing that $\nabla\muh=(\overline{\nabla\mu})_h$,
we obtain the estimates
\Beq
  \norma\muh_{\L2V}
  \leq \norma\mu_{\Lm2V}
  \leq c
  \aand
  \norma\Rh_{\L2H}
  \leq \norma R_{\Lm2H}
  \leq c \,.
  \label{stimemuhRh}
\Eeq
Now, we are ready to perform our procedure that leads to the desired regularity of the solution.
Clearly, in order to show that $\dt u\in\L2H$ and that $\mu\in\L\infty H$,
it suffices to prove estimates for the proper norms
of the different quotient $\dh u$ and of the mean~$\muh$, respectively.
To this end, we remind the reader that \eqref{primalim} and \eqref{secondalim}
have been extended up to $t=-1$.
So, we can}
integrate \eqref{primalim} with respect to time over $(s-h,s)$,
with $s\in(0,T)$,
and divide by~$h$.
We obtain
\Beq
  \iO \dh u(s) \, v
  + \iO \nabla\muh(s) \cdot \nabla v
  = \iO \Rh(s) v
  \quad \hbox{for almost every $s\in(0,T)$ and every $v\in V$}
  \non
\Eeq
and choose $v=\dh u(s)$.
Then, we integrate over $(0,t)\subset(0,T)$ with respect to~$s$ and have
\Beq
  \intQt |\dh u|^2
  + \intQt \nabla\muh \cdot \nabla\dh u
  = \intQt \Rh \, \dh u .
  \label{daprimalim}
\Eeq
At the same time, we derive from \eqref{secondalim}
\Beq
  \dh\mu
  = - \Delta\dh{\gianni u} + \dh\calW'(u)
  \quad \aeQ .
  \non
\Eeq
We multiply this equation by $\muh$ and integrate over~$Q_t$.
We have
\Beq
  \intQt \dt\muh \, \muh
  = \intQt \nabla\dh u \cdot \nabla\muh
  + \intQt \dh\calW'(u) \, \muh \,.
  \label{dasecondalim}
\Eeq
Finally, {\pier we sum \eqref{daprimalim} and \eqref{dasecondalim}}.
Two integrals cancel out and we obtain
\Beq
  \intQt |\dh u|^2
  + \frac 12 \iO |\muh(t)|^2
  = \frac 12 \iO |\muh(0)|^2
  + \intQt \Rh \, \dh u
  + \intQt \dh\calW'(u) \, \muh \,.
  \label{sommadalim}
\Eeq
Now, we estimate the integrals on the \rhs, separately.
The first one {\gianni is bounded by the first inequality in~\eqref{stimeprolmuR}.
Moreover, the second condition in \eqref{stimemuhRh} implies}
\Beq
  \intQt \Rh \, \dh u
  \leq \frac 14 \intQt |\dh u|^2
  + \intQ |\Rh|^2
  \leq \frac 14 \intQt |\dh u|^2
  + c \,.
  \non
\Eeq
For the last term of \eqref{sommadalim} some more work has to be done.
We notice that {\pier \eqref{ovunque} and} the mean value theorem yield
$\dh\calW'(u)= \calW''(\tilde u)\dh u$,
where $\tilde u$ is some {\pier measurable} function taking values
{\pier in between $u$ and of~$u(\cpto {-}h)$.}
In particular, we have $|\tilde u|\leq|u|+|u(\cpto{-}h)|$ \aeQ,
so that our assumption~\eqref{classico},
the Sobolev inequality~{\pier \eqref{sobolev}} and \eqref{stimeprol} imply for $s\in(0,T)$
\Bsist
  && \norma{\calW''(\tilde u(s))}_3^3
  \leq c \iO \bigl( 1 + |\tilde u(s)|^6 \bigr)
  \leq c \iO \bigl( 1 + |u(s)|^6 + |u(s-h)|^6 \bigr)
  \non
  \\
  && \leq c \bigl( 1 + \normaV{u(s)}^6 + \normaV{u(s-h)}^6 \bigr)
  \leq c \bigl( 1 + \norma u_{\Lm\infty V}^6 \bigr)
  = c \,.
  \non
\Esist
Hence, {\pier owing to {\gianni\eqref{stimemuhRh}} as well, we} have
\Bsist
  && \intQt \dh\calW'(u) \, \muh
  \leq \intQt |\calW''(\tilde u)| \, |\dh u| \, |\muh|
  \non
  \\
  && \leq \iot \norma{\calW''(\tilde u)}_3 \, \norma{\dh u(s)}_2 \, \norma{\muh(s)}_6 \, ds
  \leq c \iot \norma{\dh u(s)}_2 \, \normaV{\muh(s)} \, ds
  \non
  \\
  && \leq \frac 14 \iot \norma{\dh u(s)}_2^2 \, ds
  + c \ioT \normaV{\muh(s)}^2 \, ds
  \leq \frac 14 \intQt |\dh u|^2 + c \,.
  \non
\Esist
By combining all this and \eqref{sommadalim},
we conclude that the estimate
\Beq
  {\rev \intQt |\dh u|^2 + 
  \iO |\muh(t)|^2}  
  \leq c
  \non
\Eeq
holds true for $h$ small.
{\gianni
This proves that
$\dt u\in\L2H$ and $\mu\in\L\infty H$.
Now, by}
comparison in \eqref{primalim} and on account of~\eqref{defR},
we infer that $\Delta\mu\in\L2H$.
By elliptic regularity
we deduce that $\mu\in\L2W$.
Finally, by writing \eqref{secondalim} in the form
\Beq
  - \Delta u + \beta(u)
  = \mu - \lambda'(u)
  \in \L\infty H
  \non
\Eeq
and using a standard argument,
we see that $\Delta u\in\L\infty H$.
Therefore also the regularity $u\in\L\infty W$ holds true
and the proof of~\accorpa{regmubis}{regubis} is complete.
To conclude~$iii)$, we show uniqueness, i.e.,
we pick two solutions $(\mu_i,u_i,\sigma_i)$, $i=1,2$, to the limit problem
and prove that they are the same.
We set for convenience
\Beq
  \mu := \mu_1 - \mu_2 \,, \quad
  u := u_1 - u_2
  \aand
  \sigma := \sigma_1 - \sigma_2 .
  \non
\Eeq
{\pier Let us} write equations \accorpa{primalim}{terzalim} for both solutions
{\pier and exploit} the regularity result just proved
{\pier
in the course of the proof.
Taking the difference we}
obtain \aeQ
\Bsist
  && \dt u - \Delta\mu
  = R
  \label{primadiff}
  \\
  && R = \bigl( p(u_1) - p(u_2) \bigr) (\sigma_1 - \gamma\mu_1) + p(u_2) (\sigma - \gamma\mu)
  \label{diffR}
  \\
  && \mu = - \Delta u + \calW'(u_1) - \calW'(u_2)
  \label{secondadiff}
  \\
  && \dt\sigma - \Delta\sigma
  = - R
  \label{terzadiff}
\Esist
as~well as the homogeneous initial and boundary conditions.
Next, we multiply equations \eqref{primadiff}, \eqref{secondadiff} and \eqref{terzadiff}
by $u$, $\mu$ and~$\sigma$, respectively, integrate over~$Q_t$ with $t\in(0,T)$ and sum {\rev them} up.
As two integrals cancel out (thanks to boundary conditions), we~have
\Bsist
  && \frac 12 \iO |u(t)|^2
  + \intQt |\mu|^2
  + \frac 12 \iO |\sigma(t)|^2
  + \intQt |\nabla\sigma|^2
  \non
  \\
  && = \intQt R \, (u - \sigma)
  + \intQt \bigl( \calW'(u_1) - \calW'(u_2) \bigr) \mu
  \label{difftestate}
\Esist
and {\rev we} can estimate the term on the \rhs, separately, as follows.
In the sequel, the (variable) values of $c$ are allowed to depend
also on the solutions we are considering, since they are fixed.
{\pier Accounting} for~\eqref{diffR} and for the boundedness and the \Lip\ continuity of~$p$, {\pier we deduce that}
\Bsist
  && \intQt R \, (u - \sigma)
  \leq c \intQt |u| \, |\sigma_1 - \gamma\mu_1| \, (|u| + |\sigma|)
  + c \intQt |\sigma - \gamma\mu| \, (|u| + |\sigma|)
  \non
  \\
  && \leq \frac 14 \intQt |\mu|^2
  + c \iot \bigl( 1 + \norma{\sigma_1(s)-\gamma\mu_1(s)}_\infty \bigr) \bigl( \norma{u(s)}_2^2 + \norma{\sigma(s)}_2^2 \bigr) \, ds \,.
  \non
\Esist
In order to estimate the last term of~\eqref{difftestate},
we recall that $u_1$ and $u_2$ are bounded
and that $\rev \calW'$ is \Lip\ continuous on every bounded interval.
Hence, we~have
\Beq
  \intQt \bigl( \calW'(u_1) - \calW'(u_2) \bigr) \mu
  \leq c \intQt |u| \, |\mu|
  \leq \frac 14 \intQt |\mu|^2
  + c \intQt |u|^2 \,.
  \non
\Eeq
At this point, we collect \eqref{difftestate} and the above inequalities
and apply the Gronwall lemma
by noting that the function $s\mapsto\norma{\sigma_1(s){-}\gamma\mu_1(s)}_\infty$
belongs to $L^2(0,T)$
(since both $\sigma_1$ and $\mu_1$ belong to $\L2W\subset\L2\Linfty$).
This immediately leads to $u=\mu=\sigma=0$ and the proof is complete.



\vspace{3truemm}

\Begin{thebibliography}{10}

\bibitem{Barbu}
{\sc V. Barbu},
``Nonlinear semigroups and differential equations in Banach spaces'',
Noordhoff International Publishing, Leyden, 1976.

\bibitem{Brezis}
{\sc H. Brezis,}
``Op\'erateurs maximaux monotones et semi-groupes de contractions
dans les espaces de Hilbert'',
North-Holland Math. Stud.~{\bf 5},
North-Holland, Amsterdam, 1973.

{\rev
\bibitem{CH} 
{\sc J.W.~Cahn and J.E.~Hilliard}, 
{\em Free energy of a nonuniform system I. Interfacial free energy}, 
J. Chem. Phys.,
{\bf 2} (1958), pp.~258--267.
}

\bibitem{CriFriLiLowMacSanWisZheBOOK-CH2008}
{\sc V.~Cristini, H.B.~Frieboes, X.~Li, J.S.~Lowengrub, P.~Macklin, S.~Sanga,
  S.M.~Wise, and X.~Zheng}, {\em Nonlinear modeling and simulation of tumor
  growth}, in Selected Topics in Cancer Modeling: Genesis, Evolution, Immune
  Competition, and Therapy, N.~Bellomo, M.~Chaplain, and E.~{De Angelis}, eds.,
  Modeling and Simulation in Science, Engineering and Technology,
  Birkh{\"{a}}user, 2008, ch.~6, pp.~113--181.

\bibitem{CriLiLowWisJMB2009}
{\sc V.~Cristini, X.~Li, J.S.~Lowengrub, and S.M.~Wise}, {\em Nonlinear
  simulations of solid tumor growth using a mixture model: invasion and
  branching}, J.~Math. Biol., {\bf 58} (2009), pp.~723--763.

\bibitem{CriLowBOOK2010}
{\sc V.~Cristini and J.S.~Lowengrub}, ``Multiscale Modeling of Cancer: An
  Integrated Experimental and Mathematical Modeling Approach'', Cambridge
  University Press, Cambridge, 2010.

\bibitem{CriLowNieJMB2003}
{\sc V.~Cristini, J.S.~Lowengrub, and Q.~Nie}, {\em Nonlinear simulation of
  tumor growth}, J.~Math. Biol., {\bf 46} (2003), pp.~191--224.

\bibitem{DHP}
{\gianni {\sc R. Denk, M. Hieber, and J. Pr\"uss},
{\em Optimal $L^p$-$L^q$-estimates for parabolic
boundary value problems with inhomogeneous data},
Math. Z., {\bf 257} (2007), pp.~193--224.}

{\rev
\bibitem{EllSt} 
{\sc C.M. Elliott and A.M. Stuart}, 
{\em Viscous Cahn-Hilliard equation. II. Analysis}, 
J. Differential Equations,
{\bf 128} (1996),  pp.~387--414.
}

{\rev
\bibitem{EllZh} 
{\sc C.M. Elliott and S. Zheng}, 
{\em On the Cahn-Hilliard equation}, 
Arch. Rational Mech. Anal.,
{\bf 96} (1986), pp.~339--357.
}

\bibitem{FGR}
{\sc S. Frigeri, M. Grasselli, and E. Rocca},
{\em{\rev Well-posedness, regularity and global
attractor for a diffuse-interface tumor-growth model}},
in preparation.

{\rev \bibitem{HDVDZO}
{\sc A.~Hawkins-Daarud, K.G.~van der Zee, and J.T.~Oden,}, 
{\em Numerical simulation of a thermodynamically consistent four-species tumor growth model}, Int. J. Numer. Methods Biomed. Eng., {\bf 28} (2012), pp.~3--24.
}

\bibitem{HKNZ}
{\sc D.~Hilhorst, J.~Kampmann, T.N.~Nguyen, and K.G.~van der Zee}, 
{\em Formal asymptotic limit of a diffuse-interface tumor-growth model}, 
in preparation.

\bibitem{Lions}
{\sc J.-L.~Lions},
``Quelques m\'ethodes de r\'esolution des probl\`emes
aux limites non lin\'eaires'',
Dunod; Gauthier-Villars, Paris, 1969.

\bibitem{LowFriJinChuLiMacWisCriNL2010}
{\sc J.S.~Lowengrub, H.B.~Frieboes, F.~Jin, Y.-L.~Chuang, X.~Li, P.~Macklin,
  S.M.~Wise, and V.~Cristini}, {\em Nonlinear modeling of cancer: bridging the
  gap between cells and tumors}, Nonlinearity, {\bf 23} (2010), pp.~R1--R91.

  \bibitem{OdeHawPruM3AS2010}
{\sc J.T.~Oden, A.~Hawkins, and S.~Prudhomme}, {\em General diffuse-interface
  theories and an approach to predictive tumor growth modeling}, Math. Models
  Methods Appl. Sci., {\bf 20} (2010), pp.~477--517.

\bibitem{Showalter}
{\sc R.E. Showalter},
``Monotone operators in Banach spaces and nonlinear partial differential equations'', {\rev Math. Surveys Monogr.~{\bf 49},  American Mathematical Society, Providence, RI, 1997.}

\bibitem{Simon}
{\sc J. Simon},
{\em Compact sets in the space $L^p(0,T; B)$},
Ann. Mat. Pura Appl.~(4), {\bf 146} (1987), pp.~65-96.

\bibitem{Zheng}
{\sc S. Zheng}, ``Nonlinear evolution equations'',
{\pier Chapman Hall/CRC Monogr. Surv. Pure Appl. Math.~{\bf 133},
Chapman \& Hall/CRC, Boca Raton, FL, 2004.}

\End{thebibliography}

\End{document}

\bye